\documentclass{article}

\input{tcilatex}
\begin{document}

\textbf{Computation of a Feynman integral and some identities of Clausen
function values.}

\ \ \ \ \ \ \ \ \ \ \ \ \ \ \ \ \ \ \ \ \ \ \ \ \ \ \ \ \ \ \ \ \ \ \ \ \ \
\ \ \ \ \ \ \ \ \ \ \ \ \ Gert Almkvist

\textbf{Introduction.}

For the 3-loop tetrahedral Feynman diagram with non-adjacent lines carrying
masses $a$ and $b$ one gets the integral ( see Broadhurst [1] )%
\[
C(a,b)=-\frac{16}{b}\dint\limits_{2}^{2+b}\frac{dw}{w(w+a)\sqrt{w^{2}+b^{2}-4%
}}\arctan h\frac{w^{2}-4-2b}{w\sqrt{w^{2}+b^{2}-4}}
\]%
\[
-\frac{16}{b}\dint\limits_{2+b}^{\infty }\frac{dw}{w(w+a)\sqrt{w^{2}+b^{2}-4}%
}\arctan h\frac{b}{\sqrt{w^{2}+b^{2}-4}}
\]%
Broadhurst also made the conjecture%
\[
C(1,1)=4\sqrt{2}(Cl_{2}(4\alpha )-Cl_{2}(2\alpha ))
\]%
where $\sin (\alpha )=\frac{1}{3}$ and%
\[
Cl_{2}(x)=-\dint\limits_{0}^{x}\log (2\left\vert \sin (\frac{\varphi }{2}%
)\right\vert )d\varphi =\sum_{n=1}^{\infty }\frac{\sin (nx)}{n^{2}}
\]%
is the Clausen function. I computed the two integrals ( for $a=b=1$
)\bigskip by hand and found using the LLL-algorithm the following identities

\textbf{Conjecture 1. }Let \ $\tan (\alpha )=\frac{1}{\sqrt{2}}$ \ and \ $%
\tan (\beta )=\sqrt{8}+\sqrt{3}$\textbf{\ \ }Then

\textbf{1.1. }$Cl_{2}(\alpha )+Cl_{2}(\pi -\alpha )+Cl_{2}(\frac{\pi }{3}%
-\alpha )-Cl_{2}(\frac{2\pi }{3}-\alpha )=\frac{7}{4}Cl_{2}(\frac{2\pi }{3})$

\textbf{1.2.} $Cl_{2}(6\alpha -\pi )+Cl_{2}(\pi +2\alpha )-2Cl_{2}(2\alpha
)+2Cl_{2}(\pi -4\alpha )=0$

\textbf{1.3. }$Cl_{2}(\pi -2\beta )+Cl_{2}(2\beta -4\alpha )+Cl_{2}(2\beta
-2\alpha )-Cl_{2}(2\beta +2\alpha -\pi )$

$-Cl_{2}(2\alpha )-2Cl_{2}(\pi -4\alpha )-2Cl_{2}(\pi +2\alpha )=0$

\textbf{1.4. }$-12Cl_{2}(2\beta -2\alpha )+4Cl_{2}(\pi -4\alpha
)-12Cl_{2}(\pi -2\beta )-18Cl_{2}(\pi +2\alpha )+7Cl_{2}(4\alpha )=0$%
\[
\]%
The first three identities give Broadhurst's value of \ $C(1,1).$ But I
could not prove any of the identities.

In July 1998 I met Broadhurst in Vancouver and he gave me [2] where he finds
the value of \ $C(a,b)$ \ for general \ $a$ \ and \ $b$ \ and also proves it
is correct. But I wanted to prove 1.1-1.4 and decided to compute \ $C(a,b)$
\ for \ $a^{2}+b^{2}<4$ \ using elementary methods. The result contained 32
different Clausen values. Putting \ $a=\frac{1}{\pi }$ \ and \ $b=\frac{1}{e}
$ \ and using the LLL-algorithm I found two general relations that
specialized to 1.1 and 1.3 \ when \ \ $a=b=1.$ These two identities were
easily proved by differentiation. 1.2 follows from an identity in [2] \ but
1.4 \ seems very difficult to prove (it is not a consequence of any general
identity I found)%
\[
\]

2.\textbf{Computation of the integrals.}

We start with two Lemmata.

\textbf{Lemma 1.} If \ $a^{2}+b^{2}\geq c^{2}$ \ and 
\[
\delta _{1}=2\arctan (\frac{a+c}{b+\sqrt{a^{2}+b^{2}-c^{2}}})
\]%
\[
\delta _{2}=2\arctan (\frac{a+c}{-b+\sqrt{a^{2}+b^{2}-c^{2}}})
\]%
then%
\[
\dint\limits_{\alpha }^{\beta }\log (a\cos (\varphi )+b\sin (\varphi
)+c)d\varphi 
\]%
\[
=(\beta -\alpha )\log (\frac{\sqrt{a^{2}+b^{2}}}{2})+Cl_{2}(\delta
_{2}-\beta )-Cl_{2}(\delta _{2}-\alpha )+Cl_{2}(\delta _{1}+\alpha
)-Cl_{2}(\delta _{1}+\beta )
\]

\textbf{Proof: }We have%
\[
a\cos (\varphi )+b\sin (\varphi )+c=2\sqrt{a^{2}+b^{2}}\sin (\frac{\delta
_{2}-\varphi }{2})\sin (\frac{\delta _{1}+\varphi }{2}) 
\]%
Taking logarithms and integrating we are done.

\textbf{Lemma 2: }%
\[
\dint\limits_{\alpha }^{\beta }\log (\tan (\varphi )-\tan (\delta ))d\varphi 
\]%
\[
=\frac{1}{2}\left\{ Cl_{2}(2\alpha -2\delta )-Cl_{2}(2\beta -2\delta
)+Cl_{2}(\pi -2\alpha )-Cl_{2}(\pi -2\beta )\right\} -(\beta -\alpha )\log
(\cos (\delta )) 
\]

\textbf{Proof: }We have%
\[
\log (\tan (\varphi )-\tan (\delta ))=\log (2\sin (\varphi -\delta ))-\log
(2\cos (\varphi ))-\log (\cos (\delta )) 
\]%
Then use 
\[
\dint\limits_{\alpha }^{\beta }\log (2\cos (\varphi ))d\varphi =\frac{1}{2}%
\left\{ Cl_{2}(\pi -2\beta )-Cl_{2}(\pi -2\alpha )\right\} 
\]%
\[
\]

Using the partial fraction%
\[
\frac{1}{w(w+a)}=\frac{1}{a}(\frac{1}{w}-\frac{1}{w+a}) 
\]%
we get two integrals not containing \ $a$%
\[
I_{1}=\dint\limits_{2+b}^{\infty }\frac{dw}{w\sqrt{w^{2}+b^{2}-4}}\arctan h(%
\frac{b}{\sqrt{w^{2}+b^{2}-4}}) 
\]%
\[
I_{2}=\dint\limits_{2}^{2+b}\frac{dw}{w\sqrt{w^{2}+b^{2}-4}}\arctan h(\frac{%
w^{2}-2(2+b)}{w\sqrt{w^{2}+b^{2}-4}}) 
\]%
Put \ $c=\sqrt{4-b^{2}}$ \ and make the substitution%
\[
w=\frac{c}{\sin (\varphi )} 
\]%
and use the notation%
\[
\sin (\alpha _{1})=\sqrt{\frac{2-b}{2+b}} 
\]%
\[
\tan (\alpha _{2})=\frac{c}{b} 
\]%
Then%
\[
I_{1}=\frac{1}{2c}\dint\limits_{0}^{\alpha _{1}}\log (\frac{c\cos (\varphi
)+b\sin (\varphi )}{c\cos (\varphi )-b\sin (\varphi )})d\varphi 
\]%
\[
=\frac{1}{4c}\left\{ -Cl_{2}(2\alpha _{1}+2\alpha _{2})+Cl_{2}(2\alpha
_{1}-2\alpha _{2})+2Cl_{2}(2\alpha _{2})\right\} 
\]%
\[
=\frac{1}{4c}\left\{ -q_{1}+q_{2}+2q_{3}\right\} 
\]%
Similarly using Lemma 1 and 2%
\[
I_{2}=\frac{1}{2c}\dint\limits_{\alpha _{1}}^{\alpha _{2}}\left\{ \log (%
\frac{\cos (\varphi )-\frac{b}{2}}{\cos (\varphi )+\frac{b}{2}})-2\log (\tan
(\frac{\varphi }{2}))\right\} 
\]%
\[
=\frac{1}{2c}\left\{ 
\begin{array}{c}
-Cl_{2}(\alpha _{2}-\alpha _{1})+Cl_{2}(\alpha _{2}+\alpha
_{1})-Cl_{2}(2\alpha _{2})-Cl_{2}(\pi -2\alpha _{2})+Cl_{2}(\pi -\alpha
_{1}-\alpha _{2}) \\ 
-Cl_{2}(\pi +\alpha _{1}-\alpha _{2})+2Cl_{2}(\alpha _{2})-2Cl_{2}(\alpha
_{1})+2Cl_{2}(\pi -\alpha _{2})-2Cl_{2}(\pi -\alpha _{1})%
\end{array}%
\right\} 
\]%
\[
=\frac{1}{2c}\left\{
-q_{4}+q_{5}-q_{6}-q_{7}+q_{8}-q_{9}+2q_{10}-2q_{11}+2q_{12}-2q_{13}\right\} 
\]%
We will show that 
\[
I_{1}+I_{2}=0 
\]%
Using the duplication formula%
\[
Cl_{2}(2x)=2Cl_{2}(x)-2Cl_{2}(\pi -x) 
\]%
we find%
\[
q_{1}=2q_{5}-2q_{8} 
\]%
\[
q_{2}=2q_{9}-2q_{4} 
\]%
and also 
\[
q_{3}=q_{6} 
\]%
It follows that \ $I_{1}+I_{2}=0$ \ is equivalent to%
\[
2q_{4}+q_{7}-2q_{8}-2q_{10}+2q_{11}-2q_{12}+2q_{13}=0 
\]%
which follows from the following result \ ( $\alpha =\alpha _{1},\beta
=\alpha _{2}$ )

\textbf{Theorem 1: }Assume \ $\sin (\alpha )=\tan (\frac{\beta }{2})$ \ Then%
\[
Cl_{2}(\pi -2\beta )-2Cl_{2}(\beta )-2Cl_{2}(\pi -\beta )+2Cl_{2}(\alpha
)+2Cl_{2}(\pi -\alpha )+2Cl_{2}(\beta -\alpha )-2Cl_{2}(\pi -\alpha -\beta
)=0 
\]

\textbf{Proof: }Since \ $Cl_{2}(\pi )=0$ \ we find that the identity is true
when \ $\alpha =\beta =0.$ Let \ $\tan (\frac{\alpha }{2})=t$ \ and consider
\ $\alpha $ and \ $\beta $ \ as functions of \ $t.$ We have%
\[
\frac{d}{d\alpha }Cl_{2}(\alpha )=-\log (2\sin (\frac{\alpha }{2}))
\]%
\[
\frac{d}{d\alpha }Cl_{2}(\pi -\alpha )=\log (\cos (\frac{\alpha }{2}))
\]%
Let \ $LHS=f(t)$. \ Then%
\[
-2\frac{df}{dt}=\left\{ -2\log (2\cos (\beta ))-2\log (2\sin (\frac{\beta }{2%
}))+2\log (2\cos (\frac{\beta }{2}))\right\} \beta ^{\prime }
\]%
\[
+\left\{ 2\log (2\sin (\frac{\alpha }{2}))-2\log (2\cos (\frac{\alpha }{2}%
))\right\} \alpha ^{\prime }+2\log (2\sin (\frac{\beta -\alpha }{2}))(\beta
^{\prime }-\alpha ^{\prime })
\]%
\[
+2\log (2\cos (\frac{\beta +\alpha }{2}))(\beta ^{\prime }+\alpha ^{\prime })
\]%
\[
=\log (\frac{2\cos (\frac{\beta }{2})\sin (\frac{\beta -\alpha }{2})\cos (%
\frac{\beta +\alpha }{2})}{\sin (\frac{\beta }{2})\cos (\beta )})\beta
^{\prime }+\log (\frac{\sin (\frac{\alpha }{2})\cos (\frac{\beta +\alpha }{2}%
)}{\cos (\frac{\alpha }{2})\sin (\frac{\beta -\alpha }{2})})\alpha ^{\prime }
\]%
Using \ $\tan (\frac{\alpha }{2})=t$ \ we obtain%
\[
\cos (\alpha )=\frac{1-t^{2}}{1+t^{2}}
\]%
\[
\sin (\alpha )=\tan (\frac{\beta }{2})=\frac{2t}{1+t^{2}}
\]%
\[
\cos (\beta )=\frac{(1-t^{2})^{2}}{1+6t^{2}+t^{4}}
\]%
\[
\sin (\beta )=\frac{4t(1+t^{2})}{1+6t^{2}+t^{4}}
\]%
Then%
\[
\frac{2\cos (\frac{\beta }{2})\sin (\frac{\beta -\alpha }{2})\cos (\frac{%
\beta +\alpha }{2})}{\sin (\frac{\beta }{2})\cos (\beta )}=\frac{\sin (\beta
)-\sin (\alpha )}{\tan (\frac{\beta }{2})\cos (\beta )}
\]%
\[
=\frac{\frac{4t(1+t^{2})}{1+6t^{2}+t^{4}}-\frac{2t}{1+t^{2}}}{\frac{2t}{%
1+t^{2}}\cdot \frac{(1-t^{2})^{2}}{1+6t^{2}+t^{4}}}=1
\]%
Similarly%
\[
\frac{\sin (\frac{\alpha }{2})\cos (\frac{\beta +\alpha }{2})}{\cos (\frac{%
\alpha }{2})\sin (\frac{\beta -\alpha }{2})}=1
\]%
Hence \ $\frac{df}{dt}=0$ \ and since \ $f(0)=0$ $\ $we have \ $f(t)\equiv 0$%
\[
\]

\textbf{Corollarium (Conjecture 1.1): }\ Assume \ $\sin (\alpha )=\frac{1}{%
\sqrt{3}}$. Then%
\[
Cl_{2}(\alpha )+Cl_{2}(\pi -\alpha )+Cl_{2}(\frac{\pi }{3}-\alpha )-Cl_{2}(%
\frac{2\pi }{3}-\alpha )=\frac{7}{4}Cl_{2}(\frac{2\pi }{3}) 
\]

\textbf{Proof: }We get \ 
\[
\tan (\frac{\beta }{2})=\sin (\alpha )=\frac{1}{\sqrt{3}} 
\]%
and hence \ $\beta =\frac{\pi }{3}.$ Using \ $Cl_{2}(\frac{\pi }{3})=\frac{3%
}{2}Cl_{2}(\frac{2\pi }{3})$ \ and the Theorem we get the wanted result.%
\[
\]

Since \ $I_{1}+I_{2}=0$ \ we get%
\[
C(a,b)=\frac{16}{ab}\left\{ I_{3}+I_{4}\right\} 
\]%
where%
\[
I_{3}=\dint\limits_{2+b}^{\infty }\frac{dw}{(w+a)\sqrt{w^{2}-c^{2}}}\arctan
h(\frac{b}{\sqrt{w^{2}-c^{2}}})
\]%
\[
I_{4}=\dint\limits_{2}^{2+b}\frac{dw}{(w+a)\sqrt{w^{2}-c^{2}}}\arctan h(%
\frac{w^{2}-4-2b}{w\sqrt{w^{2}-c^{2}}})
\]%
In \ $I_{3}$ \ we substitute 
\[
w+a=\frac{1}{u}
\]%
Let \ 
\[
p=a+b+2
\]%
\[
d=\sqrt{4-a^{2}-b^{2}}
\]%
Then we get%
\[
I_{3}=\dint\limits_{0}^{1/p}\frac{du}{\sqrt{1-2au-d^{2}u^{2}}}\arctan h(%
\frac{bu}{\sqrt{1-2au+d^{2}u^{2}}})
\]%
\[
=\frac{1}{2d}\dint\limits_{\alpha _{3}}^{\alpha _{4}}\log (\frac{dc\cos
(\varphi )+bc\sin (\varphi )-ab}{dc\cos (\varphi )-bc\sin (\varphi )+ab}%
)d\varphi 
\]%
after the substitution%
\[
u=\frac{c\sin (\varphi )-a}{d^{2}}
\]%
where%
\[
\sin (\alpha _{3})=\frac{a}{c}
\]%
\[
\sin (\alpha _{4})=\frac{a}{c}+\frac{d^{2}}{cp}
\]%
Let%
\[
\delta _{1}=2\arctan (\frac{cd-ab}{2d+bc})
\]%
\[
\delta _{2}=2\arctan (\frac{cd-ab}{2d-bc})
\]%
\[
\delta _{3}=2\arctan (\frac{cd+ab}{2d-bc})
\]%
\[
\delta _{4}=2\arctan (\frac{cd+ab}{2d+bc})
\]%
Using Lemma 2 we obtain%
\[
I_{3}=\frac{1}{2d}\left\{ 
\begin{array}{c}
Cl_{2}(\delta _{2}-\alpha _{4})-Cl_{2}(\delta _{2}-\alpha
_{3})+Cl_{2}(\delta _{1}+\alpha _{3})-Cl_{2}(\delta _{1}+\alpha _{4}) \\ 
-Cl_{2}(\delta _{4}-\alpha _{4})+Cl_{2}(\delta _{4}-\alpha
_{3})-Cl_{2}(\delta _{3}+\alpha _{3})+Cl_{2}(\delta _{3}+\alpha _{4})%
\end{array}%
\right\} 
\]%
\[
=\frac{1}{2d}\left\{ r_{1}-r_{2}+r_{3}-r_{4}-r_{5}+r_{6}-r_{7}+r_{8}\right\} 
\]%
Finally we substitute 
\[
cu=w+\sqrt{w^{2}-c^{2}}
\]%
to get%
\[
I_{4}=\frac{1}{c}\dint\limits_{u_{1}}^{u_{2}}\frac{du}{u^{2}+\frac{2a}{c}u+1}%
\log (\frac{u^{2}-f^{2}}{f^{2}u^{2}-1})
\]%
where%
\[
f^{2}=\frac{2+b}{2-b}
\]%
\[
u_{1}=f
\]%
\[
u_{2}=f+\sqrt{\frac{2b}{2-b}}
\]%
Then put%
\[
u=\frac{d\tan (\varphi )-a}{c}
\]%
to get%
\[
I_{4}=\frac{1}{d}\dint\limits_{\alpha _{6}}^{\alpha _{7}}\log (\frac{%
d^{2}(\tan (\varphi )-\frac{a}{d})^{2}(\tan (\varphi )-\frac{a+b+2}{d})(\tan
(\varphi )-\frac{a-b-2}{d})}{c^{2}f^{2}(\tan (\varphi )-\frac{a-b+2}{d}%
)(\tan (\varphi )-\frac{a+b-2}{d})})d\varphi 
\]%
where%
\[
\tan (\alpha _{6})=\frac{p}{d}
\]%
\[
\tan (\alpha _{7})=\frac{p+\sqrt{2b^{2}+4b}}{d}
\]%
Define%
\[
\delta _{7}=\arctan (\frac{a}{d})
\]%
\[
\delta _{8}=\arctan (\frac{a+b+2}{d})=\alpha _{6}
\]%
\[
\delta _{9}=\arctan (\frac{a-b-2}{d})
\]%
\[
\delta _{10}=\arctan (\frac{a-b+2}{d})
\]%
\[
\delta _{11}=\arctan (\frac{a+b-2}{d})
\]%
Using Lemma 2 we get (all logarithms cancel)%
\[
I_{4}=\frac{1}{2d}\left\{ 
\begin{array}{c}
2Cl_{2}(2\alpha _{6}-2\delta _{7})-2Cl_{2}(2\alpha _{7}-2\delta
_{7})-Cl_{2}(2\alpha _{7}-2\delta _{8})+Cl_{2}(2\alpha _{6}-2\delta _{9}) \\ 
-Cl_{2}(2\alpha _{7}-2\delta _{9})-Cl_{2}(2\alpha _{6}-2\delta
_{10})+Cl_{2}(2\alpha _{6}-2\delta _{11})+Cl_{2}(2\alpha _{7}-2\alpha _{11})
\\ 
+2Cl_{2}(\pi -2\alpha _{6})-2Cl_{2}(\pi -2\alpha _{7})%
\end{array}%
\right\} 
\]%
\[
=\frac{1}{2d}\left\{
2r_{9}-2r_{10}-r_{11}+r_{12}-r_{13}-r_{14}+r_{15}-r_{16}+r_{17}+2r_{18}-2r_{19}\right\} 
\]%
Using LLL on $\left\{ r_{1},r_{2},...,r_{19}\right\} $ \ with the special
values \ $a=\frac{1}{\pi }$ \ and \ $b=\frac{1}{e}$ we found%
\[
r_{2}=r_{9}
\]%
\[
r_{5}=r_{11}
\]%
\[
r_{4}=-r_{13}
\]%
\[
r_{1}=r_{15}
\]%
\[
r_{8}=-r_{17}
\]%
\[
r_{6}=r_{18}
\]%
These identities are valid for all \ $a$ \ and \ $b$ \ such that \ $%
a^{2}+b^{2}<4.$ \ They are not trivial. E.g. \ $r_{5}=r_{11}$ says that%
\[
\frac{1}{2}(\delta _{4}-\alpha _{4})=\alpha _{7}-\delta _{8}
\]%
which is equivalent to (after taking $\tan $ \ on both sides)%
\[
\frac{du}{p^{2}+d^{2}+pu}=\frac{\frac{ab+cd}{2d+bc}-\frac{pc-du}{d^{2}+ap}}{%
1+\frac{(ab+cd)(pc-du)}{(2d+bc)(d^{2}+ap)}}
\]%
where%
\[
c=\sqrt{4-b^{2}}
\]%
\[
d=\sqrt{4-a^{2}-b^{2}}
\]%
\[
p=a+b+2
\]%
\[
u=\sqrt{2b^{2}+4b}
\]%
Maple verifies this identity and the other relations as well. Using these
identities we reduce \ $I_{3}+I_{4}$ \ to eleven Clausen values. In order to
simplify this and to get the same answer as Broadhurst got in [2] we have to
introduce some more notation. Let%
\[
\phi =\arctan (\frac{d}{p})
\]%
\[
\phi _{a}=\arctan (\frac{d}{a})
\]%
\[
\phi _{b}=\arctan (\frac{d}{b})
\]%
\[
s_{1}=Cl_{2}(4\phi )
\]%
\[
s_{2}=Cl_{2}(2\phi _{a}+2\phi _{b}-2\phi )
\]%
\[
s_{3}=Cl_{2}(2\phi _{a}-2\phi )
\]%
\[
s_{4}=Cl_{2}(2\phi _{b}-2\phi )
\]%
\[
s_{5}=Cl_{2}(2\phi _{a}+2\phi _{b}-4\phi )
\]%
\[
s_{6}=Cl_{2}(2\phi _{a})
\]%
\[
s_{7}=Cl_{2}(2\phi _{b})
\]%
\[
s_{8}=Cl_{2}(2\phi )
\]%
Then Broadhurst's result is 
\[
2d(I_{3}+I_{4})=s_{1}+s_{2}+s_{3}+s_{4}-s_{5}-s_{6}-s_{7}-s_{8}
\]%
Comparing with our value for \ $I_{3}+I_{4}$ \ this is equivalent to the
following identities found by LLL%
\[
-2r_{10}-2r_{11}-r_{14}+2r_{15}-2r_{19}-s_{1}+s_{6}+4s_{8}=0
\]%
\[
r_{3}=s_{4}
\]%
\[
r_{7}=-s_{2}
\]%
\[
r_{9}=s_{3}
\]%
\[
r_{12}=-s_{7}
\]%
\[
r_{16}=s_{5}
\]%
\[
r_{18}=s_{8}
\]%
The last six of these follow from%
\[
\alpha _{3}=\frac{\pi }{2}-\phi _{a}
\]%
\[
\alpha _{6}=\frac{\pi }{2}-\phi 
\]%
\[
\delta _{1}=-2\phi +\phi _{a}+2\phi _{b}-\frac{\pi }{2}
\]%
\[
\delta _{3}=2\phi -\phi _{a}-2\phi _{b}+\frac{3\pi }{2}
\]%
\[
\delta _{7}=\frac{\pi }{2}-\phi _{a}
\]%
\[
\delta _{9}=-\phi +\phi _{b}-\frac{\pi }{2}
\]%
\[
\delta _{11}=\frac{\pi }{2}+\phi -\phi _{a}-\phi _{b}
\]%
which are easily verified by Maple. It remains to prove the first identity.
It follows from the following result.

\textbf{Proposition 1:} Let \ 
\[
d=\sqrt{4-a^{2}-b^{2}} 
\]%
\[
p=a+b+2 
\]%
\[
\gamma =\arctan (\frac{p+\sqrt{2b^{2}+4b}}{d}) 
\]%
\[
\phi =\arctan (\frac{d}{p}) 
\]%
\[
\phi _{a}=\arctan (\frac{d}{a}) 
\]%
Then%
\[
2Cl_{2}(2\gamma +2\phi _{a}-\pi )+2Cl_{2}(2\gamma +2\phi -\pi )+Cl_{2}(2\phi
_{a}-4\phi )-2Cl_{2}(2\gamma -2\phi +2\phi _{a}-\pi ) 
\]%
\[
+2Cl_{2}(\pi -2\gamma )+Cl_{2}(4\phi )-Cl_{2}(2\phi _{a})-4Cl_{2}(2\phi )=0 
\]

\textbf{Proof:} consider the left hand side as a function \ $f(a,b).$
Differentiating we get%
\[
-2df=4\log (2\sin (\gamma +\phi _{a}-\frac{\pi }{2}))(d\gamma +d\phi
_{a})+4\log (2\sin (\gamma +\phi -\frac{\pi }{2}))(d\gamma +d\phi ) 
\]%
\[
2\log (2\sin (\phi _{a}-2\phi ))(d\phi _{a}-2d\phi )-4\log (2\sin (\gamma
-\phi +\phi _{a}-\frac{\pi }{2}))(d\gamma -d\phi +d\phi _{a}) 
\]%
\[
-4\log (2\sin (\frac{\pi }{2}-\gamma ))d\gamma +4\log (2\sin (2\phi ))d\phi
-2\log (2\sin (\phi _{a}))d\phi _{a}-8\log (2\sin (\phi ))d\phi 
\]%
\[
=4\log (\frac{\sin (\gamma +\phi _{a}-\frac{\pi }{2})\sin (\gamma +\phi -%
\frac{\pi }{2})}{\sin ^{2}(\gamma -\phi +\phi _{a}-\frac{\pi }{2})\sin (\phi
_{a})})d\gamma 
\]%
\[
+2\log (\frac{\sin ^{2}(\gamma +\phi _{a}-\frac{\pi }{2})\sin (\phi
_{a}-2\phi )}{\sin ^{2}(\gamma -\phi +\phi _{a}-\frac{\pi }{2})\sin (\phi
_{a})})d\phi _{a} 
\]%
\[
+4\log (\frac{\sin (\gamma +\phi -\frac{\pi }{2})\sin (\gamma -\phi +\phi
_{a}-\frac{\pi }{2})\sin (2\phi )}{\sin (\phi _{a}-2\phi )\sin ^{2}(\phi )}%
)d\phi 
\]%
\[
=4\log (T_{1})d\gamma +2\log (T_{2})d\phi _{a}+4\log (T_{3})d\phi =0 
\]%
since \ $T_{1}=T_{2}=T_{3}=1.$ We show that \ $T_{3}=1$ which is equivalent
to%
\[
\cos (\gamma +\phi )\cos (\gamma +\phi _{a}-\phi )\sin (2\phi )=\sin (\phi
_{a}-2\phi )\sin ^{2}(\phi ) 
\]%
if and only if%
\[
2(1-\tan (\gamma )\tan (\phi ))(1-\tan (\gamma )\tan (\phi _{a})+\tan (\phi
)\tan (\gamma )+\tan (\phi )\tan (\phi _{a})) 
\]%
\[
=\tan (\phi )(1+\tan ^{2}(\gamma ))(\tan (\phi _{a})(1-\tan ^{2}(\phi
))-2\tan (\phi )) 
\]%
if and only if%
\[
2p\sqrt{2b^{2}+4b}(p^{2}-d^{2}-2ap+(b+2)\sqrt{2b^{2}+4b}) 
\]%
\[
=(p^{2}-d^{2}-2ap)(d^{2}+p^{2}+2b^{2}+4b+2p\sqrt{2b^{2}+4b}) 
\]%
which is easily verified. Similarly one verifies that \ $T_{1}=T_{2}=1$. \
Hence $df=0$ which means that \ $f(a,b)$ \ is constant. Putting \ $a=b=0$ \
we get \ $d=2,$ $p=2,$ $\gamma =\phi =\frac{\pi }{4},$ $\phi _{a}=\frac{\pi 
}{2}$ \ which gives \ $f(0,0)=0.$ \ It follows that $\ f(a,b)\equiv 0$ \ and
the proof is finished.

\textbf{Corollary:} We have if \ $a^{2}+b^{2}<4$%
\[
C(a,b)=\frac{8}{ab\sqrt{4-a^{2}-b^{2}}}\left\{ 
\begin{array}{c}
Cl_{2}(4\phi )+Cl_{2}(2\phi _{a}+2\phi _{b}-2\phi )+Cl_{2}(2\phi _{a}-2\phi
)+Cl_{2}(2\phi _{b}-2\phi ) \\ 
-Cl_{2}(2\phi _{a}+2\phi _{b}-4\phi )-Cl_{2}(2\phi _{a})-Cl_{2}(2\phi
_{b})-Cl_{2}(2\phi )%
\end{array}%
\right\} 
\]

\textbf{Proposition 2: }We have%
\[
2Cl_{2}(2\phi )-4Cl_{2}(2\phi _{b})+Cl_{2}(4\phi _{b})+2Cl_{2}(2\phi
_{b}-2\phi )-2Cl_{2}(2\phi _{a}-2\phi ) 
\]%
\[
+Cl_{2}(2\phi _{a}-4\phi )+2Cl_{2}(2\phi _{a}+2\phi _{b}-2\phi
)-Cl_{2}(2\phi _{a}+4\phi _{b}-4\phi )=0 
\]

\textbf{Proof: }Let \ $f(a,b)$ \ be the LHS. Differentiating we get%
\[
df=4\log (\frac{\sin (\phi )\sin (\phi _{a}-\phi )\sin (\phi _{a}+2\phi
_{b}-2\phi )}{\sin (\phi _{b}-\phi )\sin (\phi _{a}-2\phi )\sin (\phi
_{a}+\phi _{b}-\phi )})d\phi 
\]%
\[
+2\log (\frac{\sin (\phi _{a}-2\phi )\sin ^{2}(\phi _{a}+\phi _{b}-\phi )}{%
\sin ^{2}(\phi _{b}-\phi )\sin (\phi _{a}+2\phi _{b}-2\phi )})d\phi _{a}
\]%
\[
+4\log (\frac{\sin (2\phi _{b})\sin (\phi _{b}-\phi )\sin (\phi _{a}+\phi
_{b}-\phi )}{\sin ^{2}(\phi _{b})\sin (\phi _{a}+2\phi _{b}-2\phi )})d\phi
_{b}=0
\]%
since all logarithms are zero (checked by Maple). Since \ $f(0,0)=0$ \ we
have \ $f(a,b)\equiv 0.$

\textbf{Corollarium 1: }Identity 1.2 is true.

\textbf{Proof: }Put \ $a=b=1$ \ in Proposition 2. Let \ $\alpha =\arctan (%
\frac{1}{\sqrt{2}})$. \ Then \ $\phi =\frac{\pi }{2}-2\alpha $ \ and \ $\phi
_{a}=\phi _{b}=\frac{\pi }{2}-\alpha $ . Inserting this into Proposition 2
we obtain%
\[
2CL_{2}(\pi -4\alpha )-4Cl_{2}(\pi -2\alpha )-Cl_{2}(4\alpha
)+Cl_{2}(6\alpha +\pi )-Cl_{2}(\pi +2\alpha )=0 
\]%
Using the duplication formula%
\[
Cl_{2}(4\alpha )=2Cl_{2}(2\alpha )-2Cl_{2}(\pi -2\alpha ) 
\]%
we get 1.2

\textbf{Corollarium 2. }Identity 1.3 is true

\textbf{Proof: \ }\ With \ $a=b=1$ \ in Proposition 1 we have \ $\tan
(\gamma )=\sqrt{8}+\sqrt{3}$ .\ Hence \ $\gamma =\beta $ \ in 1.3.
Proposition 1 gives%
\[
2Cl_{2}(2\beta -2\alpha )+2Cl_{2}(2\beta -4\alpha )+Cl_{2}(6\alpha -\pi
)-2Cl_{2}(2\beta +2\alpha -\pi )+2Cl_{2}(\pi -2\beta ) 
\]%
\[
+Cl_{2}(2\pi -8\alpha )-Cl_{2}(\pi -2\alpha )-4Cl_{2}(\pi -4\alpha )=0 
\]%
Using 1.2 for \ $Cl_{2}(6\alpha -\pi )$ \ and the duplication formula for \ $%
Cl_{2}(8\alpha )$ \ we obtain 1.3

\textbf{Acknowledgements. }The beginning of this work was done in 1998 at
CECM, Simon Frazer University. I want to thank Jonathan and Peter Borwein
for many stimulating discussions. The final part of the paper could not have
been written without Broadhurst's result in [2]-

\textbf{References:}

\textbf{1. }D.Broadhurst, A dilogarithmic 3-dimensional Ising tetrahedron,
hep-th/9805025

\textbf{2.} D.Broadhurst, Solving differential equations for 3-loop
diagrams: relation to hyperbolic geometry and knot theory, hep-th/9806174.

\bigskip

\textbf{Appendix. Proof of an identity found by Broadhurst}

Let \ $\sin (\alpha )=\frac{1}{3}$. \ Then Broadhurst in [1] found the
following identity using PSLQ%
\[
4\sqrt{2}(Cl_{2}(4\alpha )-Cl_{2}(2\alpha )) 
\]%
\[
=\sum_{n=0}^{\infty }(-\frac{1}{2})^{3n}\frac{1}{n+\frac{1}{2}}(\frac{1}{n+%
\frac{1}{2}}-3\log (2))-3\sum_{n=1}^{\infty }(-\frac{1}{2})^{3n}\frac{H_{n}}{%
n+\frac{1}{2}} 
\]%
where%
\[
H_{n}=\sum_{k=1}^{n}\frac{1}{k} 
\]%
We will prove this here using several formulas for the dilogarithm%
\[
Li_{2}(x)=\sum_{n=1}^{\infty }\frac{x^{n}}{n^{2}} 
\]%
We start by quoting some formulas in Lewin [2] pp.244

\textbf{Lemma 1}

\textbf{1.1.}%
\[
Li_{2}(x)+Li_{2}(-x)=\frac{1}{2}Li_{2}(x^{2}) 
\]

\textbf{1.2.}%
\[
Li_{2}(x)+Li_{2}(-\frac{x}{1-x})=-\frac{1}{2}\log ^{2}(1-x) 
\]

\textbf{1.3.}%
\[
Li_{2}(\frac{1}{1+x})-Li_{2}(-x)=\frac{\pi ^{2}}{6}-\frac{1}{2}\log
(1+x)\log (\frac{(1+x)}{x^{2}}) 
\]

\textbf{1.4.}(Abel)%
\[
Li_{2}(\frac{x}{1-x}\cdot \frac{y}{1-y})=Li_{2}(\frac{x}{1-y})+Li_{2}(\frac{y%
}{1-x})-Li_{2}(x)-Li_{2}(y)-\log (1-x)\log (1-y) 
\]

\textbf{1.5.}%
\[
Li_{2}(x)+Li_{2}(1-x)=\frac{\pi ^{2}}{6}-\log (x)\log (1-x) 
\]

\textbf{Lemma 2. }We have 
\[
\sum_{n=1}^{\infty }\frac{H_{n}}{2n+1}z^{2n+1} 
\]%
\[
=\frac{1}{2}\left\{ \frac{1}{2}\log ^{2}(1-z)-\frac{1}{2}\log ^{2}(1+z)+\log
(2)\log (\frac{1-z}{1+z})+Li_{2}(\frac{1+z}{2})-Li_{2}(\frac{1-z}{2}%
)\right\} 
\]

\textbf{Proof: }We have by [3], p.715%
\[
\sum_{n=1}^{\infty }H_{n}x^{2n}=-\frac{\log (1-x^{2})}{1-x^{2}} 
\]%
Integrate%
\[
\sum_{N01}^{\infty }\frac{H_{n}}{2n+1}x^{2n+1}=-\frac{1}{2}%
\dint\limits_{0}^{x}(\log (1-t)+\log (1+t))(\frac{1}{1-t}+\frac{1}{1+t})dt 
\]%
and we are done by formula 3.2 p.266 in Lewin [2]%
\[
\dint\limits_{0}^{x}\frac{\log (a+bt)}{c+et}dt 
\]%
\[
=\frac{1}{e}\left\{ \log (\frac{ae-bc}{e}\log (\frac{c+ex}{c})-Li_{2}(\frac{%
b(c+ex)}{bc-ae})+Li_{2}(\frac{bc}{bc-ae}\right\} 
\]

\textbf{Proof of Broadhurst's identity.}

Introduce the notation%
\[
x=\frac{1}{2}(1+\frac{i}{\sqrt{8}}) 
\]%
\[
y=\frac{1}{2} 
\]%
\[
u=\frac{\sqrt{8}+i}{3} 
\]%
\[
z=-\frac{i}{8} 
\]%
Then we want to compute 
\[
Cl_{2}(4\alpha )-Cl_{2}(2\alpha )=\func{Im}(Li_{2}(u^{4})-Li_{2}(u^{2})) 
\]%
In Abel's identity (Lemma 1.4) we have%
\[
\frac{x}{1-x}=u^{2},\text{ }\frac{y}{1-y}=1,\text{ }\frac{x}{1-y}=1-z,\text{ 
}\frac{y}{1-x}=\frac{1}{1+z} 
\]%
and hence

(2.1) \ $Li_{2}(u^{2})=Li_{2}(1-z)+Li_{2}(\frac{1}{1+z})-Li_{2}(x)-Li_{2}(%
\frac{1}{2})+\log (2)\log (1-x)$\ \ \bigskip

By Lemma 1.5 and 1.3 we get

(2.2) $\ Li_{2}(1-z)=-Li_{2}(z)+\frac{\pi ^{2}}{6}-\frac{1}{2}\log (z)\log
(1-z)$

(2.3) \ $Li_{2}(\frac{1}{1+z})=Li_{2}(-z)+\frac{\pi ^{2}}{6}-\frac{1}{2}\log
(1+z)\log (\frac{1+z}{z^{2}})$

Adding (2.1), (2.2) and (2.3) we obtain

(2.4) \ $Li_{2}(u^{2})=Li_{2}(\overline{z})-Li_{2}(z)-Li_{2}(x)+\frac{\pi
^{2}}{3}-Li_{2}(\frac{1}{2})+\log (2)\log (1-x)$

$-\log (z)\log (1-z)-\frac{1}{2}\log (1+z)\log (\frac{1+z}{z^{2}})$

By Lemma 1.2 and 1.1

(2.5) \ $Li_{2}(x)+Li_{2}(-u^{2})=-\frac{1}{2}\log ^{2}(1-x)$

(2.6) \ $Li_{2}(u^{2})+Li_{2}(-u^{2})=\frac{1}{2}Li_{2}(u^{4})$

(2.7) \ $Li_{2}(u^{4})=2Li_{2}(u^{2})-2Li_{2}(x)-\log ^{2}(1-x)$

Adding (2.4) and (2.7) gives

(2.8) \ $Li_{2}(u^{4})-Li_{2}(u^{2})=Li_{2}(\overline{z}%
)-Li_{2}(z)-3Li_{2}(x)+\frac{\pi ^{2}}{3}-Li_{2}(\frac{1}{2})+\log (2)\log
(1-x)-\log ^{2}(1-x)$

\ \ \ \ \ \ \ \ \ \ \ \ \ \ \ \ \ \ \ \ \ \ \ \ \ \ \ \ \ \ \ \ \ \ \ \ $%
-\log (z)\log (1-z)-\frac{1}{2}\log (1+z)\log (\frac{1+z}{z^{2}})$

We want

(2.9) \ $\func{Im}(Li_{2}(u^{4})-Li_{2}(u^{2}))=\frac{1}{i}(Li_{2}(\overline{%
z})-Li_{2}(z))-\frac{3}{2i}(Li_{2}(x)-Li_{2}(\overline{x}))$

\ \ \ \ \ \ \ \ \ \ \ \ \ \ \ \ \ \ \ \ \ \ \ \ \ \ \ \ \ \ \ \ \ \ \ \ \ \
\ \ \ \ \ \ $+\func{Im}\left\{ \log (2)\log (1-x)-\log ^{2}(1-x)-\log
(z)\log (1-z)-\frac{1}{2}\log (1+z)\log (\frac{1+z}{z^{2}})\right\} $

We get%
\[
\frac{1}{i}(Li_{2}(\overline{z})-Li_{2}(z))=\frac{1}{i}\sum_{n=1}^{\infty }%
\frac{1}{n^{2}}\left\{ (\frac{i}{\sqrt{8}})^{n}-(\frac{-i}{\sqrt{8}}%
)^{n}\right\} =\frac{1}{\sqrt{2}}\sum_{n=0}^{\infty }\frac{1}{(2n+1)^{2}}(-%
\frac{1}{2})^{3n} 
\]%
and using Lemma 2%
\[
\frac{1}{2i}(Li_{2}(x)-Li_{2}(\overline{x}))=\frac{1}{2i}(Li_{2}(\frac{1-z}{2%
})-Li_{2}(\frac{1+z}{2})) 
\]%
\[
=-\frac{z}{i}\sum_{n=1}^{\infty }\frac{H_{n}}{2n+1}z^{2n}+\frac{1}{4i}%
\left\{ Li_{2}(\frac{1-z}{2})-Li_{2}(\frac{1+z}{2})\right\} +\frac{\log (2)}{%
2i}\log (\frac{1-z}{1+z}) 
\]%
\[
=\frac{1}{\sqrt{8}}\sum_{n=1}^{\infty }\frac{H_{n}}{2n+1}(-\frac{1}{2}%
)^{3n}+\alpha \log (\sqrt{\frac{9}{8}})+\frac{\log (2)}{\sqrt{8}}%
\sum_{n=0}^{\infty }\frac{1}{2n+1}(-\frac{1}{2})^{3n} 
\]%
Putting everything together we obtain Broadhurst's formula (all logarithms
cancel).

\textbf{References.}

\textbf{1. }D.Broadhurst, A dilogarithmic 3-dimensional tetrahedron,
hep-th/980525

\textbf{2. }L.Lewin, Dilogarithms and associated functions, Macdonald,
London 1958\bigskip

\textbf{3. }A.P.Prudnikov, O.I.Marichev, Yu.A.Brychkov, Integrals and series
I (Russian) Moskva Nauka, 1986

\bigskip

Institute of Algebraic Meditation

Fogdar\"{o}d 208

S-2433 H\"{o}\"{o}r

Sweden

gert@maths.lth.se

\end{document}